\documentclass[11pt,english]{article}
\usepackage{graphicx} 
\usepackage{amssymb, amsthm, amsmath, fullpage, hyperref}
\usepackage{graphicx, enumerate}
\usepackage{tikz-network}
\usetikzlibrary{positioning, 
                quotes}
\usepackage[utf8]{inputenc}
\usepackage{cleveref}
\crefname{thm}{Theorem}{Theorems}
\crefname{thm}{Theorem}{Theorems}
\crefname{lem}{Lemma}{Lemmas}
\usepackage{amsfonts}
\usepackage{amssymb}
\usepackage{MnSymbol}
\usepackage{amsmath,amsthm,dsfont}
\usepackage[english]{babel}
\usepackage{makeidx}
\usepackage{mathtools}
\usepackage{bbm}
\usepackage{fullpage,amsfonts,amssymb,epsfig,epstopdf,amsmath,titling,url,array,titlesec}
\usepackage{mathtools}
\usepackage{authblk}
\usepackage{tikz}
\usepackage{float}
\usepackage{makeidx}
\usepackage[mathscr]{euscript}
\usepackage{nomencl}
\usepackage{comment}
\theoremstyle{plain}
\newtheorem{thm}{Theorem}[section]
\newtheorem{lem}[thm]{Lemma}
\newtheorem{cor}[thm]{Corollary}

\theoremstyle{definition}
\newtheorem{defn}[thm]{Definition}

\newtheorem{prop}[thm]{Proposition}
\newtheorem{rem}[thm]{Remark}

\newtheorem{obs}[thm]{Observation}
\newtheorem{claim}[thm]{Claim}
\newtheorem*{notn}{Notation}

\newtheorem{case}{Case}
\makeatletter
\@addtoreset{case}{subsection} 
 
\makeatother
\providecommand{\keywords}[1]
{
  \textbf{Keywords:} #1
}

\author[1]{Rajat Adak}
\author[2]{L. Sunil Chandran}
\affil[1]{\texttt{rajatadak@iisc.ac.in}} \affil[2]{\texttt{sunil@iisc.ac.in}}
\affil[1,2]{Department of Computer Science and Automation,
Indian Institute of Science, Bangalore, India}
\date{}

\title{Vertex-Based Localization of generalized Tur\'{a}n Problems}

\begin{document}

\maketitle
\begin{abstract}
Let $\mathcal{F}$ be a family of graphs. A graph is called $\mathcal{F}$-free if it does not contain any member of $\mathcal{F}$. Generalized Tur\'{a}n problems aim to maximize the number of copies of a graph $H$ in an $n$-vertex $\mathcal{F}$-free graph. This maximum is denoted by $ex(n, H, \mathcal{F})$. When $H \cong K_2$, it is simply denoted by $ex(n,F)$. Erdős and Gallai established the bounds $ex(n, P_{k+1}) \leq \frac{n(k-1)}{2}$ and $ex(n, C_{\geq k+1}) \leq \frac{k(n-1)}{2}$. This was later extended by Luo \cite{luo2018maximum}, who showed that $ex(n, K_s, P_{k+1}) \leq \frac{n}{k} \binom{k}{s}$ and $ex(n, K_s, C_{\geq k+1}) \leq \frac{n-1}{k-1} \binom{k}{s}$. Let $N(G,K_s)$ denote the number of copies of $K_s$ in $G$. In this paper, we use the vertex-based localization framework, introduced in \cite{adak2025vertex}, to generalize Luo's bounds. In a graph $G$, for each $v \in V(G)$, define $p(v)$ to be the length of the longest path that contains $v$. We show that 
\[N(G,K_s) \leq \sum_{v \in V(G)} \frac{1}{p(v)+1}{p(v)+1\choose s} = \frac{1}{s}\sum_{v \in V(G)}{p(v) \choose s-1}\]
We strengthen the cycle bound from \cite{luo2018maximum} as follows: In graph $G$, for each $v \in V(G)$, let $c(v)$ be the length of the longest cycle that contains $v$, or $2$ if $v$ is not part of any cycle. We prove that  
\[N(G,K_s) \leq \left(\sum_{v\in V(G)}\frac{1}{c(v)-1}{c(v) \choose s}\right) - \frac{1}{c(u)-1}{c(u) \choose s}\]  
where $c(u)$ denotes the circumference of $G$.  Furthermore, we characterize the class of extremal graphs that attain equality for these bounds. We provide full proofs for the cases $s = 1$ and $s \geq 3$, while the case $s = 2$ follows from the result in \cite{adak2025vertex}. We also conclude with a generalization of a result by Balister–Bollobás–Riordan–Schelp \cite{BALISTER2003366}.

\end{abstract}

\section{Introduction}
In recent years, there has been growing interest in generalized Turán problems, which focus on determining the maximum number of copies of a fixed graph $H$ that a graph on $n$ vertices can contain under various constraints. This line of inquiry can be traced back to the classical theorem of Tur\'{a}n~\cite{Turan}, where $H$ is assumed to be an edge, in which cases, the goal is just maximizing the number of edges in a $n$-vertex graph under the given constraints; the maximum value in this case is referred to as the Tur\'{a}n number of $n$ for the given constraint.

\subsection{Classical Tur\'{a}n problems}
A graph is said to be $\mathcal{F}$-free if it contains no members of the family $\mathcal{F}$ as a subgraph. The Turán number of such a graph on $n$ vertices is denoted by $ex(n, T)$. A classical example of a Turán problem was given by Erd\H{o}s and Gallai~\cite{gallai1959maximal}, who studied the Turán numbers of graphs that avoid long paths and long cycles — that is, $ex(n, P_{k+1})$ and $ex(n, C_{\geq k+1})$, where $P_{k+1}$ is the path on $k+1$ vertices and $C_{\geq k+1}$ denotes any cycle of length at least $k+1$.

\begin{thm}\label{thm:path}\emph{(Erd\H{o}s-Gallai \cite{gallai1959maximal})}
For a simple graph $G$ with $n$ vertices without a path of length $k (\geq 1)$, that is a path with $k$ edges $(P_{k+1})$,
\begin{equation*}
    |E(G)| \leq ex(n,P_{k+1})\leq \frac{n(k-1)}{2}
\end{equation*}
and equality holds if and only if all the components of $G$ are complete graphs of order $k$.
\end{thm}

\begin{thm}\label{th:clas}\emph{(Erd\H{o}s-Gallai \cite{gallai1959maximal})}
    For a simple graph $G$ with $n$ vertices and circumference (length of the longest cycle in the graph) at most $k$,
    \[|E(G)| \leq ex(n,C_{\geq k+1}) \leq \frac{k(n-1)}{2}\]
    and equality holds if and only if $G$ is connected and all its blocks are complete graphs of order $k$.
\end{thm}
\subsection{Generalized Tur\'{a}n problems}
The maximum number of copies of a graph $H$ in an $n$-vertex, $\mathcal{F}$-free graph is denoted by $ex(n, H, \mathcal{F})$. Luo~\cite{luo2018maximum} extended \Cref{thm:path,th:clas} by providing a generalized Turán-type version of the Erd\H{o}s–Gallai theorems. Recently, Chakraborti and Chen~\cite{CHAKRABORTI2024103955} classified the extremal graphs for these bounds.
\begin{notn}
    Let $N(G, K_s)$ denote the number of subgraphs of $G$ isomorphic to $K_s$.
\end{notn}
\begin{thm}\emph{(Luo \cite{luo2018maximum})}\label{luopath}
    For a simple graph $G$ with $n$ vertices without a path of length $k (\geq 1)$,
    \begin{equation*}
        N(G,K_s) \leq ex(n,K_s,P_{k+1}) \leq \frac{n}{k}{k\choose s}
    \end{equation*}
    and equality holds if and only if all the components of $G$ are complete graphs of order $k$.
\end{thm}
\begin{rem}
    Note that taking $s =2$, we get back \cref{thm:path}.
\end{rem}
\begin{thm}\emph{(Luo \cite{luo2018maximum})}\label{luocycle}
 For a simple graph $G$ with $n$ vertices and circumference at most $k$,
    \[N(G,K_s) \leq ex(n,K_s,C_{\geq k+1}) \leq \frac{n-1}{k-1}{k\choose s}\]
    and equality holds if and only if $G$ is connected, and all its blocks are complete graphs of order $k$.    
\end{thm}
\begin{rem}
    Note that if we take $s= 2$, we get back \Cref{th:clas}.
\end{rem}
\subsection{Localization}
Recently, Brada\v{c}~\cite{bradac} and Malec and Tompkins~\cite{DBLP:journals/ejc/MalecT23} introduced the concept of \emph{localization} as a means to strengthen and generalize classical extremal graph bounds. The classical extremal bounds rely on global parameters; for example, the length of the longest path in \Cref{thm:path}, or the circumference in \Cref{th:clas}.

The localization framework begins by introducing a \emph{localization function} that captures global graph parameters in a locally defined manner. For instance, while the clique number of a graph $G$ is a global property, we can localize it by defining a function $w : E(G) \rightarrow \mathbb{Z}$, where $w(e)$ denotes the size of the largest clique containing the edge $e$. Localization replaces global constraints with local weight functions, defined over specific graph parameters such as edges or vertices. Thus far, most of the literature in this area has focused on edge-based weights, where a function $w: E(G) \rightarrow \mathbb{Z}$ assigns to each edge $e$ a value determined by the structure it is part of. For example, to localize \Cref{thm:path,th:clas}, weight $w(e)$ is defined, which represents the length of the longest path or cycle containing the edge $e$, respectively. This localized approach refines classical extremal bounds by replacing global parameters with edge-specific weights.

In~\cite{bradac,DBLP:journals/ejc/MalecT23}, the authors introduced an edge-based localization of a popular version of the classical Tur\'{a}n's theorem. Several other localization results were also provided in~\cite{DBLP:journals/ejc/MalecT23}, including an edge-based variant of \Cref{thm:path}. Subsequently, Zhao and Zhang~\cite{zhao2025localized} established an edge-based localization of \Cref{th:clas}.

Recently, Kirsch and Nir~\cite{Kirsch_2024} extended the edge-based localization of Tur\'{a}n's theorem from~\cite{bradac,DBLP:journals/ejc/MalecT23} by localizing a generalization of the classical result originally proposed by Zykov~\cite{Zykov1949-qr}. Arag\~{a}o and Souza~\cite{aragao2024localised} further advanced the theory by developing a localized version of the Graph Maclaurin Inequalities~\cite{Khadzhiivanov1977-xg}.

The study of edge-based localization has also been extended to hypergraphs. In~\cite{zhao2025localized}, the authors presented an edge-based localization of the hypergraph Erd\H{o}s--Gallai theorem for paths. A spectral version of localization was explored in~\cite{liu2024variantsspectralturantheorems}, where the authors provided an edge-based localization of the spectral Tur\'{a}n inequality.

It was observed that the inequalities derived from edge-based localization, while not always intuitively appealing, often yield distribution-type statements, from which interesting results can be derived for specialized cases. These bounds are often similar in
spirit to the generalization of Sperner’s Lemma to LYM inequality, from the theory of set systems, (see chapter 3, \cite{bollobás1986combinatorics}), or like the famous Kraft inequality. As an alternative, \emph{vertex-based localization} was introduced in~\cite{adak2025vertex}, where weight functions are assigned to vertices instead of edges. This vertex-based approach often leads to more direct and conceptually appealing generalizations of classical results.

\subsection{Vertex-based Localization}
Adak and Chandran~\cite{adak2025vertex} introduced the concept of \emph{vertex-based localization}. They provided vertex-based localizations of \Cref{th:clas,thm:path}, offering a more natural and intuitive strengthening of these classical results compared to the edge-based approaches in~\cite{DBLP:journals/ejc/MalecT23} and~\cite{zhao2025localized}. In their framework, weights are assigned to the vertices of the graph as follows:

\begin{notn}
For a graph $G$, let $p(v)$ denote the weight of a vertex $v \in V(G)$, defined as:
\begin{equation*}
    p(v) = \text{length of the longest path in $G$ that contains $v$}.
\end{equation*}
\end{notn}
\begin{notn}
For a graph $G$, let $c(v)$ denote the weight of a vertex $v \in V(G)$, defined as:
\begin{equation*}
    c(v) = \text{length of the longest cycle in $G$ that contains $v$}.
\end{equation*}
If $v$ is not contained in any cycle, then set $c(v) = 2$.
\end{notn}

\begin{thm}\label{th:pathl}\emph{(Adak-Chandran \cite{adak2025vertex})}
For a simple graph $G$, 
\begin{equation*}
    |E(G)| \leq \sum_{v\in V(G)}\frac{p(v)}{2}
\end{equation*}
Equality holds if and only if every connected component of $G$ is a clique.
\end{thm}
\begin{thm}\label{th:cyclel}\emph{(Adak-Chandran \cite{adak2025vertex})}
   For a simple graph $G$,
\[|E(G)| \leq \left(\sum_{v \in V(G)} \frac{c(v)}{2}\right) - \frac{c(u)}{2}\]
where $c(u)$ is the circumference of $G$. Equality holds if and only if $G$ is a parent-dominated block graph (defined in the next section).
\end{thm} 
\begin{rem}
    It is easy to see that one can get back \Cref{thm:path,th:clas} from \Cref{th:pathl,th:cyclel} respectively.
\end{rem}
In~\cite{adak2025turan}, the same authors introduced a vertex-based localization of a well-known variant of Tur\'{a}n's classical theorem. In a subsequent work~\cite{adak2025localizationframeworkgeneralizeextremal}, they provided a vertex-based localization of Wood's result~\cite{wood} on the maximum number of cliques in graphs with bounded maximum degree.

\section{Our Results}
In this paper, we extend \Cref{th:pathl,th:cyclel} by providing bounds on the number of copies of $K_s$, where $s\geq 3$, in a graph, rather than just $K_2$ (i.e., edges), which were the focus in those earlier results. Specifically, we present vertex-based localizations of \Cref{luopath,luocycle}, thereby initiating the study of vertex-based localization of generalized Tur\'{a}n-type problems.
\begin{defn}
A \emph{block} in a graph is a maximal connected subgraph that contains no cut-vertices. A \emph{block graph} is a connected graph in which every block is a clique.
\end{defn}

\begin{figure}[H]
    \centering
    \includegraphics[width=6cm]{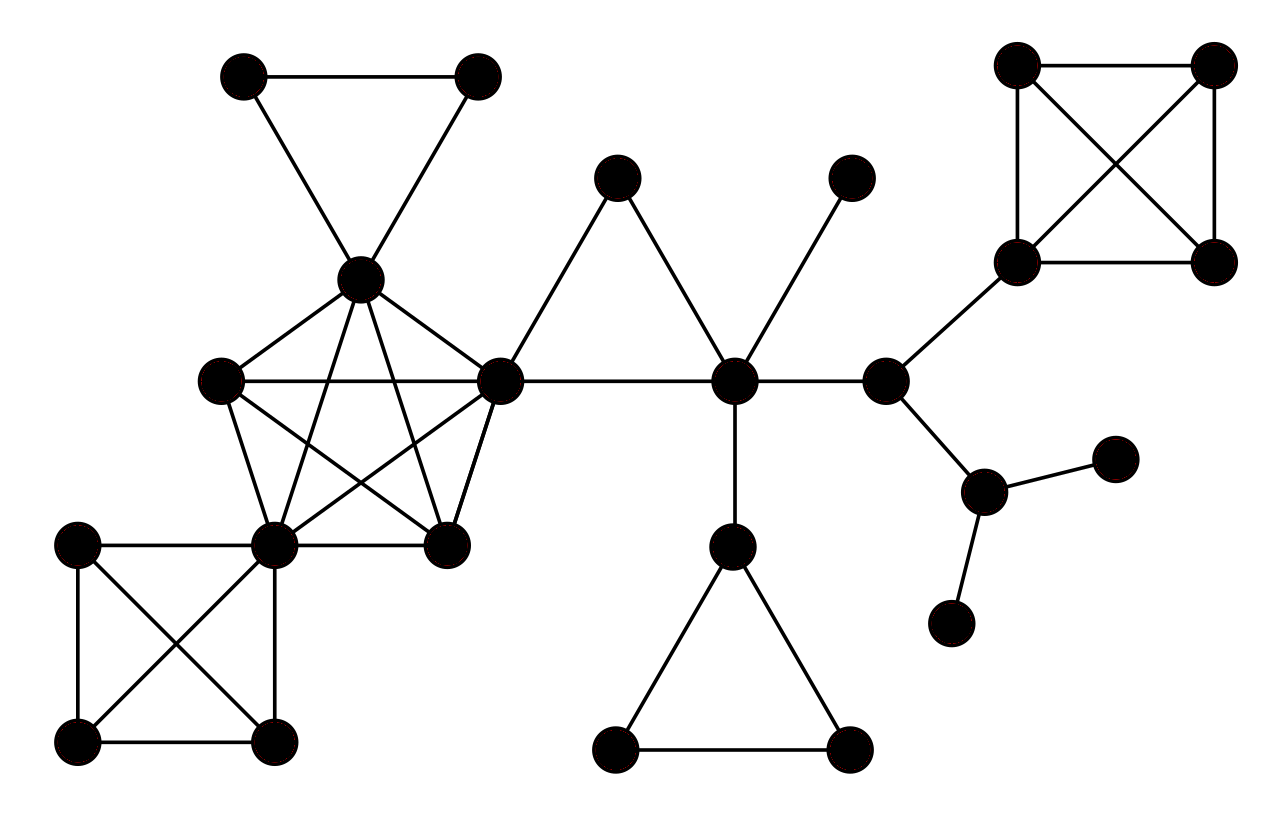}
    \caption{An example of a block graph}
    \label{fig:Path_P}
\end{figure}

\begin{rem}
A block graph $G$ can be associated with a tree $T$, known as the \emph{block tree}, where each vertex in $T$ corresponds to a block in $G$. Two vertices in $T$ are adjacent if and only if their corresponding blocks in $G$ share a vertex.

\noindent Once the block tree $T$ is defined, any vertex of $T$ can be chosen as the root, naturally inducing a parent-child relation on the blocks of $G$.
\end{rem}

\begin{defn}
A block graph $G$ is called a \emph{parent-dominated block graph} if its block tree is rooted at a vertex corresponding to a block of the largest order, and every other block $B$ in $G$ has order less than or equal to that of its parent block in the tree.
\end{defn}
\begin{defn}\label{Hc}
    Let $H_c \subseteq V(G)$ be defined as the set of vertices $v$ in $G$ for which $c(v) \geq s$, that is, 
    \[
    H_c = \{v \in V(G) \mid c(v) \geq s\}.
    \]
\end{defn}
\begin{defn}\label{Hp}
    Let $H_p \subseteq V(G)$ be the set of vertices $v$ in $G$ such that $p(v) \geq s - 1$, that is,
    \[
    H_p = \{v \in V(G) \mid p(v) \geq s - 1\}.
    \]
\end{defn}

\noindent We are now all set to formally state our main results.

\begin{thm}\label{thm1}
    Let $G$ be a simple graph with $c(u)$ as circumference. Then,
    \begin{equation*}
        N(G,K_s) \leq \left(\sum_{v\in V(G)}\frac{1}{c(v)-1}{c(v) \choose s}\right) - \frac{1}{c(u)-1}{c(u) \choose s}
    \end{equation*}
     If $s =1$, then equality holds if and only if $G$ is Hamiltonian. Otherwise, equality holds if and only if $G[H_c]$ is a Parent-dominated block graph.
\end{thm}

\begin{thm}\label{thm2}
    Let $G$ be a simple graph. Then,
    \begin{equation*}
        N(G,K_s) \leq \sum_{v \in V(G)} \frac{1}{p(v)+1}{p(v)+1\choose s} = \frac{1}{s}\sum_{v \in V(G)}{p(v) \choose s-1}
    \end{equation*}
    If $s=1$, then equality always holds. Otherwise, equality holds if and only if all the components of $G[H_p]$ are cliques.
\end{thm}
\begin{rem}
It is easy to verify that by setting $s = 2$ in \Cref{thm1,thm2}, \Cref{th:pathl,th:cyclel} can be recovered, respectively.
\end{rem}

\subsection{Proof of \Cref{thm1}}

For the proof of the inequality in \Cref{thm1}, we adopt an approach similar to that of the proof of \Cref{th:cyclel}, relying on the notions of \textit{transforms} and \textit{simple transforms}. In contrast, for the classification of the extremal graph class, we employ a simpler method. Interestingly, however, using this alternative approach, one can not classify the extremal graph class for \Cref{th:cyclel}.

\subsubsection{Some Notation, Definitions and Lemmas}

\begin{notn}
Let $P = v_0v_1\dots v_k$ be a path in a graph $G$. We define:
\[
    P(v_i,v_j) = v_iv_{i+1}\dots v_{j-1}v_j \quad \text{for } 0\leq i \leq j \leq k,
\]
\[
    P(v_i,v_j) = v_iv_{i-1}\dots v_{j+1}v_j \quad \text{for } 0\leq j \leq i \leq k.
\]
\end{notn}

\noindent A $v_0$-path in a graph $G$ is a path that begins at the vertex $v_0 \in V(G)$. Unless otherwise specified, we assume that a $v_0$-path $P$ in $G$ proceeds from $v_0$ to its other endpoint, which we refer to as the \textit{terminal vertex} of $P$. This convention will assist in defining the following functions.

We define the distance (with respect to $P$) between two vertices $x, y \in V(P)$ as the number of edges between them along the path $P$, and denote it by $dist_P(x,y)$. Note that $dist_P(x,x) = 0$.

Let $i \in \mathbb{N}$ and $s \in V(P)$ such that $dist_P(v_0, s) \geq i$. Define $Pred_s(i)$ to be the vertex $z \in V(P)$ such that $dist_P(z,s) = i$ and $v_0$ is closer to $z$ than $s$, that is, $dist_P(v_0,z) \leq dist_P(v_0,s)$. That is, $z$ is a \textit{predecessor} of $s$ on $P$ at distance $i$.

Similarly, let $j \in \mathbb{N}$ and $t \in V(P)$ such that $dist_P(t,v_k) \geq j$. Define $Succ_x(j)$ to be the vertex $z \in V(P)$ such that $dist_P(t,z) = j$ and $v_k$ is closer to $z$ than $s$, that is, $dist_P(t,v_k) \geq dist_P(z,v_k)$. That is, $z$ is a \textit{successor} of $t$ on $P$ at distance $j$.

Let $P = v_0v_1v_2\dots v_k$ be a longest $v_0$-path in $G$. It is easy to observe that $N(v_k) \subseteq \{v_0, v_1, \dots, v_{k-1}\}$; otherwise, $P$ could be extended to a longer $v_0$-path, contradicting maximality. Suppose $v_k$ is adjacent to $v_j$ for some $j \in \{0, 1, \dots, k-2\}$. Then we can construct a new path $P'$ defined as:
\[
P' = v_0v_1 \dots v_jv_kv_{k-1} \dots v_{j+2}v_{j+1},
\]
by removing the edge $v_jv_{j+1}$ and adding the edge $v_jv_k$. Note that $V(P') = V(P)$.

\begin{defn}
Let $P$ be a longest $v_0$-path in $G$, where $P = v_0v_1v_2\dots v_k$ and $v_k$ is adjacent to $v_j$ for some $0\leq j \leq k-2$. Then the path
\[
P' = v_0v_1 \dots v_jv_kv_{k-1} \dots v_{j+2}v_{j+1}
\]
is called a \textit{simple transform} of $P$.
\begin{figure}[H]
    \centering
    \includegraphics[width=12cm]{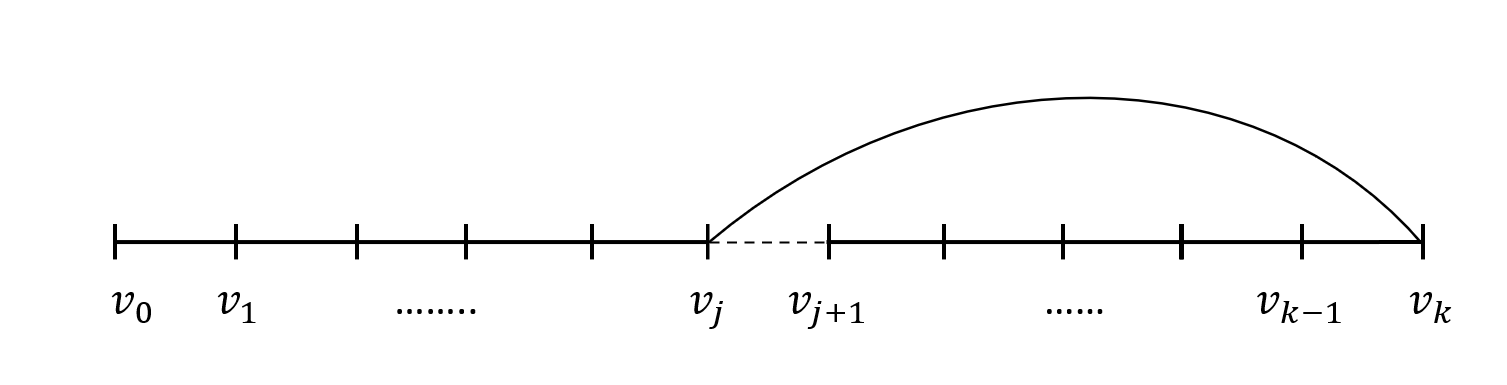}
    \caption{Path $P'$ is a \textit{simple transform} of $P$}
\end{figure}
\end{defn}

\begin{rem}\label{st_inv}
If $P'$ is a simple transform of $P$, then $P$ is a simple transform of $P'$ also.
\end{rem}

\begin{defn}
Let $P$ be a longest $v_0$-path in $G$ with terminal vertex $v_k$. Since each vertex $x \in N(v_k)$ lies on $P$, for every $x \in N(v_k) \setminus \{Pred_{v_k}(1)\}$, there exists a simple transform of $P$ ending at $Succ_x(1)$. Let $P'$ be one such simple transform. We may apply further simple transforms to $P'$, yielding additional paths. Let $P''$ be a path obtained by applying a sequence of simple transforms starting from $P$. Clearly, $P''$ is also a longest $v_0$-path in $G$, and $V(P'') = V(P)$. We refer to such a path $P''$ as a \textit{transform} of $P$. Let $\mathcal{T}$ denote the set of all transforms of $P$. Note that $\mathcal{T}$ depends on $(G, P, v_0)$. Clearly, all the paths in $\mathcal{T}$ are $v_0$-paths.
\end{defn}

\begin{rem}
If $P'$ is a transform of $P$, then $P$ is also a transform of $P'$.
\end{rem}
\begin{notn}In a graph $G$ with $P$ as a longest $v_0$-path and $\mathcal{T} = \mathcal{T}(G,P,v_0)$, we define the following notations and symbols, to be used throughout the paper.
\begin{itemize}
    \item  Let $L(G, P, v_0)$ be the set of terminal vertices of the paths in $\mathcal{T}$; that is: $$L(G,P,v_0) = \{v \in V(P) \mid v \text{ is a terminal vertex of some $v_0$-path } P' \in \mathcal{T}\}$$
     \item For $v \in L(G,P,v_0)$, let $S_v$ be the set of neighbors of $v$ outside the set $L(G,P,v_0)$, that is, $S_v = N(v) \cap (L(G,P,v_0))^c$.
      \item For $v \in L(G,P,v_0)$, let $\mathcal{T}_v$ denote the set of all the transforms in $\mathcal{T}$ with $v$ as the terminal vertex.
     \item Let $w(P_v) = Pred_v(c(v)-1)$ on $P_v \in \mathcal{T}_v$, that is, the $c(v)^{th}$ vertex on $P_v$ starting from the vertex $v$ in the first position. We call this the \textit{pivot} vertex of the path $P_v$.
     \item Define $Front^*(P_v) = P_v(v_0,w(P_v))$ and $Back^*(P_v) = P_v(w(P_v),v)$. Let $Front(P_v)$ and $Back(P_v)$ be defined from $Front^*(P_v)$ and $Back^*(P_v)$ respectively by removing the pivot vertex.
\end{itemize}
\end{notn}
\begin{rem}\label{big_cycle}
   Note that for any $v \in L(G, P, v_0)$ and $P_v \in \mathcal{T}_v$, the vertex $v$ cannot be adjacent to any vertex in $\text{Front}(P_v)$; otherwise, this would yield a cycle longer than $c(v)$ that includes $v$. Additionally, observe that the pivot vertex, $w(P_v)$, may not be adjacent to $v$.
\end{rem}
\begin{defn}
    Let $x \in L(G,P,v_0)$ such that $c(x) = min\{c(v)\mid v \in L(G,P,v_0)\}$. For $v \in L(G,P,v_0)$, a path $P_v \in \mathcal{T}_v$ is called a good path if the terminal vertex $v$ is not adjacent to any vertex in $P_v(v_0,Pred_vc(x)
    )$; otherwise $P_v$ is a bad path.
\end{defn}
\begin{lem}\label{No_early_nbr}
    For all $v \in L(G,P,v_0)$, all the paths in $\mathcal{T}_v$ are good.
    \begin{proof}
      From \cref{big_cycle}, $x$ can not have a neighbor in $Front(P_x)$. Therefore, $P_x$ is a good path. Suppose for $v \in L(G,P,v_0)$, there exists a bad path $P_v \in \mathcal{T}_v$.

        By the definition of transforms, we know that $P_v$ can be obtained from $P_x$ by applying a sequence of simple transforms. Without loss of generality, assume that $P_v$ is the first bad path in this sequence. 

        Observe that in every path in the sequence from $P_x$ to $P_v$, only the last $c(x) - 1$ vertices are being permuted, while the positions of the remaining vertices remain fixed. More precisely, let $i = |V(P)| - c(x) + 1$. Then, for every path in the sequence from $P_x$ to $P_v$, the first $i$ vertices (starting from $v_0$) are identical. The set of the last $c(x) - 1$ vertices remains the same across all these paths; only their ordering changes.

        Therefore, the set of vertices in the subpath $P_v(z_v, v)$ equals the set $\text{Back}(P_x)$, where $z_v = \text{Pred}_v(c(x) - 2)$ on path $P_v \in \mathcal{T}_v$.

        Since $P_v$ is a bad path, there exists a vertex $t \in P_v(v_0,Pred_vc(x))$, adjacent to $v$. Thus, the cycle formed by adding the edge $(v,t)$ to path $P_v(t,v)$, is of length more than $c(x)$ and contains the vertex $x$, which is a contradiction.
    \end{proof}
    \end{lem}
    
\begin{cor}\label{L in back}
    For all $v \in L(G,P,v_0)$ and $P_v \in \mathcal{T}_v$, $L(G,P,v_0) \cap V(Front^*(P_v)) = \emptyset$
        \begin{proof} Suppose there exists $v \in L(G,P,v_0)$ such that; \[L(G,P,v_0) \cap V(Front^*(P_v)) \neq \emptyset \]  Let $t \in L(G,P,v_0) \cap (V(Front^*(P_v))$. Let $P_y \in \mathcal{T}_y$ be such that one can obtain $P_t \in \mathcal{T}_t$ as a simple transform of $P_y$.  Therefore, on the path $P_y$, vertex $v$ is adjacent to $Pred_t(1)$ as shown in \cref{fig3}.
        \begin{figure}[H]
            \centering
            \includegraphics[width=16cm]{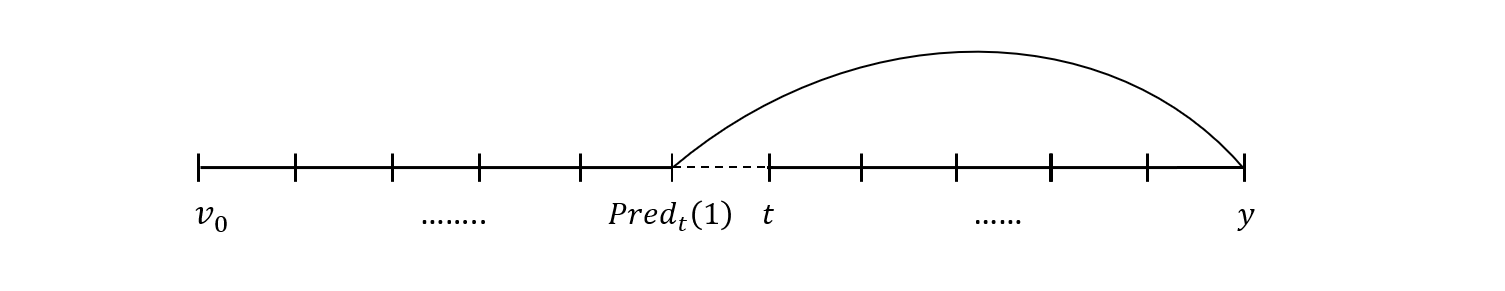}
            \caption{$P_t$ is a simple transform of $P_y$}
            \label{fig3}
        \end{figure}
        From \cref{No_early_nbr}, we got that the position of the first $i=|V(P)|-c(x)+1$ vertices on all the transforms of $P$ are invariant. Therefore, $V(Front^*(P_v)) = V(Front^*(P_y))$. Since $t \in V(Front^*(P_v))$, we get that, $Pred_t(1) \in V(Front(P_y))$. This contradicts \cref{big_cycle}, as $y$ can not be adjacent to any vertex in $Front(P_y)$.
        \end{proof}
\end{cor}
\begin{lem}
    For $v \in L(G,P,v_0)$, $|S_v| \leq c(v) - |L(G,P,v_0)|$
    \begin{proof}
         Recall that $S_v$ is the set of all the neighbors of the vertex $v$ which lie outside the set $L(G,P,v_0)$. From \cref{big_cycle}, we get that, $N(v) \subseteq Back^*(P_v)$. Thus $S_v \subseteq Back^*(P_v) \setminus L(G,P,v_0)$. From \cref{L in back} we get that $L(G,P,v_0) \cap V(Front^*(P_v)) = \emptyset$, therefore, $L(G,P,v_0) \subseteq Back(P_v) \subset Back^*(P_v)$. Therefore, we get that;
         \begin{equation}\label{s<c-l}
             |S_v| \leq |Back^*(P_v)| - |L(G,P,v_0)| = c(v) - |L(G,P,v_0)|
         \end{equation}
    \end{proof}
\end{lem}
\begin{lem}
    For $v \in L(G,P,v_0)$, $d(v) \leq |L(G,P,v_0)|$
    \begin{proof}
        Let $t \in N(v)$. Clearly $t \in V(P_v)$. Note that on path $P_v$, $Succ_t(1) \in L(G,P,v_0)$, by the definition of transforms and $L(G,P,v_0)$. Therefore, for every neighbor of $v$, there exists a unique vertex in $L(G,P,v_0)$. Thus we get;
        \begin{equation}\label{dv<l}
            d(v) \leq |L(G,P,v_0)|
        \end{equation}
    \end{proof}
\end{lem}
\subsubsection{Preprocessing for Proof of Inequality}
\begin{rem}\label{s=1}
    When $s = 1$, the quantity $N(G, K_s)$ simply counts the number of vertices in $G$, that is, $N(G, K_s) = n$. In this case, the right-hand side of the bound in \Cref{thm1} simplifies to:
    \begin{equation*}
        \left(\sum_{v \in V(G)} \frac{1}{c(v) - 1} {c(v) \choose s}\right) - \frac{1}{c(u) - 1} {c(u) \choose s} = \sum_{v \in V(G) \setminus \{u\}} \frac{c(v)}{c(v) - 1} \geq \sum_{v \in V(G) \setminus \{u\}} \frac{n}{n - 1} = n.
    \end{equation*}
    Hence, the desired bound is satisfied, and equality holds if and only if $c(v) = n$ for every $v \in V(G)$, that is, when $G$ is Hamiltonian.
\end{rem}

\begin{rem}
    When $s = 2$, observe that the statement of \cref{thm1} reduces to that of \cref{th:cyclel}. This is because $N(G, K_2) = m$, the number of edges in $G$. In this case, the bound given in \cref{thm1} becomes:
    \[
        m \leq \sum_{v \in V(G) \setminus \{u\}} \frac{c(v)}{2}.
    \]
    From \cref{Hc}, recall that $H_c = \{v \in V(G) \mid c(v) \geq s\}$. Since $c(v) \geq 2$ for all $v \in V(G)$ by definition, it follows that $H_c = V(G)$. Therefore, $G[H_c] = G$, and the extremal graph class for \cref{th:cyclel} coincides with that of \cref{thm1}. Hence, for $s = 2$, the proof of \cref{thm1} follows directly from the proof of \cref{th:cyclel} in \cite{adak2025vertex}.
\end{rem}

\noindent In what follows, we will assume that $s \geq 3$.
\subsubsection{Proof of inequality}
Let $G$ be a simple graph and let $u$ be a vertex such that $c(u) = max\{c(v) \mid v \in V(G)\}$.
\begin{notn}
   Let $N(G, K_s, X)$ denote the number of copies of $K_s$ in $G$ that contain at least one vertex from the set $X \subseteq V(G)$.
\end{notn}
\begin{prop}\label{main_claim} Let $L = L(G,P,u)$, then
    \begin{equation*}
        N(G,K_s,L) \leq \sum_{v \in L}\frac{1}{c(v)-1}{c(v) \choose s}
    \end{equation*}
\end{prop}
\noindent The proof of the above proposition will be presented after establishing a few preliminary lemmas.

\vspace{2mm}
  
   For each vertex $v \in L = L(G, P, u)$, let $N_v(G, K_s, L)$ denote the \textit{contribution} of $v$ toward the total number of $K_s$-cliques counted in $N(G, K_s, L)$. The natural question is: how should we define these individual contributions so that they add up exactly to the total count? That is, we want the following identity to hold:
\begin{equation}\label{contribution}
    N(G, K_s, L) = \sum_{v \in L} N_v(G, K_s, L)
\end{equation}

A naive approach would be to let $N_v(G, K_s, L) = N(G, K_s, \{v\})$, meaning that every copy of $K_s$ containing vertex $v$ contributes a full count of 1 to $N_v(G,K_s,L)$. However, this leads to overcounting: any clique that contains multiple vertices from $L$ will be counted once for each of those vertices. 

For example, suppose a $K_s$ subgraph $K \subseteq G$ intersects $L$ in $k$ vertices, that is, $V(K) \cap L = \{v_1, v_2, \dots, v_k\}$. Then, under the naive definition, this single clique $K$ would be counted $k$ times in the sum. 

To correct for this overcounting, we distribute the contribution of such a clique equally among the vertices it contains from $L$. That is, each vertex $v_i \in V(K) \cap L$ contributes exactly $\frac{1}{k}$ toward $N(G, K_s, L)$. 

Formally, $N_v(G, K_s, L)$ is defined as the sum over all $K_s$-cliques that contain $v$, where each such clique contributes $\frac{1}{k}$ if it contains $k$ vertices from $L$. This ensures that each clique is counted exactly once in the total sum.

\begin{rem}
    It is easy to see that, keeping all other conditions fixed, adding an extra edge to the induced subgraph $G[N(v)]$ does not decrease the value of $N_v(G, K_s, L)$. In other words, the contribution of a vertex $v$ to the total $K_s$ count is non-decreasing under edge additions within its open neighborhood. Therefore, $N_v(G, K_s, L)$ is maximized when the closed neighborhood $N[v]$ induces a clique in $G$.
    
    Vertex $v$ has $d(v)$ neighbors in total, among which $|S_v|$ lie outside the set $L(G, P, u)$. The remaining $d(v) - |S_v|$ neighbors belong to $L(G, P, u)$.

    Consider a copy $K$ of $K_s$ that includes $v$, with exactly $t$ vertices from $S_v$ and the remaining $s - t - 1$ (other than $v$) vertices from $L(G, P, u)$. Since $K$ contains $s - t$ vertices from $L(G, P, u)$, it contributes $\frac{1}{s - t}$ to $N_v(G, K_s, L)$. The number of such cliques is at most:
$${|S_v| \choose t} {d(v) - |S_v| \choose s - t - 1}.$$

Therefore, we obtain the following upper bound:
\begin{equation}\label{indi_contri}
    N_v(G, K_s, L) \leq \sum_{t = 0}^{s - 1} \frac{1}{s - t} {|S_v| \choose t} {d(v) - |S_v| \choose s - t - 1}
\end{equation}

\end{rem}

\begin{obs}\label{enough1}
    To prove \cref{main_claim}, from \cref{contribution,indi_contri} it is enough to show that for all $v \in L(G,P,u)$
    \begin{equation}\label{main_obs}
        \sum_{t=0}^{s-1}\frac{1}{s-t}{|S_v|\choose t}{d(v)-|S_v| \choose s-t-1} \leq \frac{1}{c(v)-1}{c(v) \choose s}
    \end{equation}
\end{obs}

\begin{rem}\label{|S_v|=0}
Recall from \Cref{big_cycle} that the vertex $v$ cannot be adjacent to any vertex in $Front(P_v)$, where $P_v \in \mathcal{T}_v$. Hence, the neighbors of $v$ must lie entirely within $Back^*(P_v)$, implying that $d(v) \leq |Back^*(P_v)| - 1 = c(v) - 1$. 

\noindent Now, consider the case when $|S_v| = 0$. In this case, on the right-hand side in \cref{indi_contri}, only the term that corresponds to $t =0$ remains. Thus we obtain;
\begin{equation}\label{sv =0}
    N_v(G, K_s, L) \leq \frac{1}{s} \binom{d(v)}{s-1} \leq \frac{1}{s} \binom{c(v)-1}{s-1} =\frac{1}{c(v)} \binom{c(v)}{s} \leq \frac{1}{c(v) - 1} \binom{c(v)}{s}
\end{equation}
This establishes the validity of \cref{main_claim} for the case $|S_v| = 0$. Therefore, for the remainder of the analysis, we assume $|S_v| \geq 1$.

\end{rem}

\noindent Note that we can simplify the left-hand side of \cref{main_obs} using the convolution identity as follows;
\begin{align}
    \sum_{t=0}^{s-1}\frac{1}{s-t}{|S_v|\choose t}{d(v)-|S_v| \choose s-t-1}&= \sum_{t=0}^{s-1}\frac{1}{(d(v)-|S_v|+1)}{|S_v|\choose t}{d(v)-|S_v|+1 \choose s-t}\nonumber\\
    &=\frac{1}{d(v)-|S_v|+1}\left(\left(\sum_{t=0}^{s}{|S_v|\choose t}{d(v)-|S_v|+1 \choose s-t}\right) - {|S_v|\choose s}\right)\nonumber\\
    &=\frac{1}{d(v)-|S_v|+1}\left({d(v)+1 \choose s} - {|S_v| \choose s}\right)\label{convulation}
\end{align}
Thus using \cref{convulation}, the bound in \cref{main_obs} becomes equivalent to
\begin{equation}\label{convulation1}
    \frac{1}{d(v)-|S_v|+1}\left({d(v)+1 \choose s} - {|S_v| \choose s}\right) \leq \frac{1}{c(v)-1}{c(v) \choose s}
\end{equation}
Observe that the left-hand side of \cref{convulation1} depends on $d(v)$ and $|S_v|$, while the right-hand side depends on $c(v)$. As a result, the two sides are not directly comparable in their current form. To overcome this, we aim to establish a stronger inequality that implies \cref{convulation1} and ensures that both sides are functions of the same parameter, making them directly comparable.
\begin{obs}\label{inc}
    For $s \geq 2$, it is easy to  verify that $\frac{1}{x-1}{x\choose s}$ is an non-decreasing for integral values of $x$ such that $x >1$ 
\end{obs}
 From \cref{s<c-l,dv<l} we get that $c(v) \geq |S_v| + |L(G,P,u)| \geq |S_v| + d(v)$. Since $|S_v| \geq 1$, we have $d(v) \geq 1$. Therefore $|S_v| + d(v) >1$. Thus from \cref{inc} we get that;
\begin{equation}\label{s+d<c}
    \frac{1}{|S_v|+d(v)-1}{|S_v| + d(v) \choose s}\leq \frac{1}{c(v)-1}{c(v)\choose s}
\end{equation}

\begin{obs}\label{obs26}
    In the view of \cref{s+d<c}, to establish \cref{convulation1} it is enough to show that;
    \begin{equation*}
        \frac{1}{d(v) - |S_v| + 1} \left( {d(v) + 1 \choose s} - {|S_v| \choose s} \right) \leq \frac{1}{{|S_v| + d(v) - 1}} { |S_v| + d(v) \choose s }
    \end{equation*}
\end{obs}

\begin{lem}\label{usefull}
    For $x\geq 4$ and $y \geq 3$
    \begin{equation*}
        \frac{1}{x-3}{x-2\choose y} \leq \frac{1}{x-1}{x-1\choose y}
    \end{equation*}
    Equality holds if and only if $x-1 <y$.
    \begin{proof}
        First assuming $y =2$, we get that;
        \begin{equation*}
            \frac{1}{x-3}{x-2\choose 2} = \frac{x-2}{2} = \frac{1}{x-1}\frac{(x-1)(x-2)}{2}=\frac{1}{x-1}{x-1\choose 2}
        \end{equation*}
        Now assume $y \geq 3$. If $x-2<y$, then the bound holds trivially. It is easy to check that equality holds if and only if $x-1<y$.

       \noindent Suppose $x-2\geq y$, given $y \geq 3$. For  contradiction, assume;
        \begin{align*}
            &\frac{1}{x-3}{x-2\choose y} \geq \frac{1}{x-1}{x-2\choose s}\\
            \implies &\frac{1}{x-3}\frac{(x-2)!}{y!(x-2-y)!} \geq \frac{1}{x-1}\frac{(x-1)!}{y!(x-1-y)!}\\
            &\implies (x-1-y)\geq (x-3)
        \end{align*}
        This is contradiction, since $y \geq 3 \implies (x-1-y)\leq (x-4)<(x-3)$
    \end{proof}
\end{lem}
In the following lemma, we establish the inequality from \cref{obs26} for $c \geq 3$ and $|S_v| \geq 1$. The remaining cases have already been settled.
\begin{lem}\label{enough2}
    For $v \in L(G,P,u)$, such that $d(v) \geq |S_v| \geq 1$ and with $s \geq 3$ 
    \begin{equation}\label{enough2eqn}
        \frac{1}{d(v) - |S_v| + 1} \left( {d(v) + 1 \choose s} - {|S_v| \choose s} \right) \leq \frac{1}{{|S_v| + d(v) - 1}} { |S_v| + d(v) \choose s }
    \end{equation}

\begin{proof}
Note that the bound in \cref{enough2} holds with equality when $s = 2$, and hence does not yield any \emph{non-trivial} inequality. 
This suffices to establish the inequality of \cref{thm1} for $s = 2$, but the case does not contribute towards the characterization of the extremal graphs. 
Therefore, in what follows, we restrict attention to $s \geq 3$ and refer the reader to \cite{adak2025vertex} for a complete treatment of the $s = 2$ case.  

For $s \geq 3$, the analysis produces inequalities that are not only sufficient to establish the bound in \cref{thm1}, but also play a crucial role in ruling out certain cases during the identification of the extremal class. 

 The proof proceeds by induction on the quantity $d(v) + |S_v| + s$.
    
    \vspace{1mm}
    \noindent \textit{Base Case:}  $|S_v| + d(v)+s = 5$ with $|S_v| = 1, d(v) = 1, s=3$
    \begin{equation*}
        LHS=\frac{1}{d(v) - |S_v| + 1} \left( {d(v) + 1 \choose s} - {|S_v| \choose s} \right)  = \frac{1}{1-1+1}\left({2\choose 3} - {1\choose 3}\right) = 0
    \end{equation*}
    \begin{equation*}
        RHS=\frac{1}{{|S_v| + d(v) - 1}} { |S_v| + d(v) \choose s } = \frac{1}{1+1-1}{2\choose 3} = 0
    \end{equation*}
        
    \vspace{1mm}
    \noindent\textit{Induction Step:} We will proceed by considering the following cases.
    \begin{case} $|S_v|=1$. The inequality in \cref{enough2} simplifies to
\begin{equation*}
    \frac{1}{d(v)} \binom{d(v)+1}{s} \leq \frac{1}{d(v)} \binom{d(v)+1}{s},
\end{equation*}
which clearly holds with equality.
\end{case}
    \begin{case}  Assume $|S_v| >1 \implies d(v) >1$. Note that we have;
    \begin{align}\label{LHS}
        LHS &= \frac{1}{d(v) - |S_v| + 1} \left( \binom{d(v) + 1}{s} - \binom{|S_v|}{s} \right)\nonumber\\ 
        &= \frac{1}{d(v) - |S_v| + 1} \left( {d(v)  \choose s} - {|S_v|-1 \choose s}+{d(v) \choose s-1} - {|S_v|-1 \choose s-1}\right )
    \end{align}
    \begin{align}\label{RHS}
        RHS &= \frac{1}{{|S_v| + d(v) - 1}} { |S_v| + d(v) \choose s }\nonumber \\
        &=  \frac{1}{{|S_v| + d(v) - 1}} \left({ |S_v| + d(v)-1 \choose s } + { |S_v| + d(v)-1 \choose s -1}\right)
    \end{align}
    Moreover $(d(v) - 1) \geq 1$ and $(|S_v| - 1) \geq 1$ and $s\geq 3$ as assumed. Therefore, we can apply \textit{induction hypothesis} with parameters $(d(v)-1)$ and $(|S_v| -1)$. Therefore, we have:

    \begin{equation}\label{eqnn}
        \frac{1}{(d(v)-1)- (|S_v|-1) + 1} \left( {d(v) \choose s} - {|S_v|-1 \choose s}  \right) \leq \frac{1}{{|S_v| + d(v) - 3}} { |S_v| + d(v)-2 \choose s }
    \end{equation}
    Since $d(v) + |S_v| \geq 4$ and $s \geq 3$, substituting $x = d(v) + |S_v|$ and $y = s$ in \cref{usefull}, we get;
    \begin{equation}\label{y=s}
        \frac{1}{|S_v|+d(v)-3}{|S_v|+d(v)-2 \choose s} \leq \frac{1}{(|S_v| + d(v) -1)}{|S_v| + d(v)-1 \choose s}
    \end{equation}
    and equality holds if and only if $(x-1) = |S_v| + d(v) -1 < s$.
    Therefore, from \cref{eqnn,y=s};
    \begin{equation}\label{y=s ineq}
        \frac{1}{(d(v)-1)- (|S_v|-1) + 1} \left( {d(v) \choose s} - {|S_v|-1 \choose s}  \right) \leq \frac{1}{(|S_v| + d(v) -1)}{|S_v| + d(v)-1 \choose s}
    \end{equation}
\begin{obs}\label{y=s equality}
    If the bound in \cref{y=s ineq} is tight, then $|S_v|+d(v) -1 <s$.
\end{obs}

\vspace{2mm}
\noindent Note that if $s-1 \geq 3$, using \textit{induction hypothesis} we get that;
\begin{equation}\label{s-1>3}
    \frac{1}{(d(v)-1)- (|S_v|-1) + 1} \left( {d(v) \choose s-1} - {|S_v|-1 \choose s-1}  \right) \leq \frac{1}{{|S_v| + d(v) - 3}} { |S_v| + d(v)-2 \choose s-1 }
\end{equation}
Since $|S_v| + d(v) \geq 4$ and if we assume $s-1 \geq 3$, then substituting $x = d(v) + |S_v|$ and $y =s-1$ in \cref{usefull}, we get that;
\begin{equation}\label{s-1}
    \frac{1}{(|S_v|+d(v) - 3)} {|S_v| + d(v) -2\choose s-1} \leq \frac{1}{(|S_v|+d(v) - 1)}{|S_v| + d(v) -1\choose s-1}
\end{equation}
Thus, from \cref{s-1>3,s-1}, we get that if $s-1 \geq 3$;
\begin{equation}\label{eqn19}
    \frac{1}{(d(v)- |S_v| + 1)} \left( {d(v) \choose s-1} - {|S_v|-1 \choose s-1}  \right) \leq \frac{1}{(|S_v|+d(v) - 1)}{|S_v| + d(v) -1\choose s-1}
\end{equation}
\begin{obs}\label{y=s-1 equality}
    If the bound in \cref{eqn19} is tight, then $|S_v|+d(v) -1 <s-1$.
\end{obs}

\vspace{2mm}

\noindent Now suppose $s-1 =2$. Since $d(v) \geq |S_v| \geq 2$ the left-hand side of \cref{enough2eqn} becomes; 
\begin{align*}
    \frac{1}{(|S_v|+d(v) - 1)}{|S_v| + d(v) -1\choose s-1} &= \frac{|S_v|+d(v)-2}{2}
\end{align*}
Now, if $s-1 =2$, the right-hand side of \cref{enough2eqn} simplifies to;
\begin{align*}
     \frac{1}{(d(v)- |S_v| + 1)} \left( {d(v) \choose s-1} - {|S_v|-1 \choose s-1}  \right) &= \frac{1}{(d(v)- |S_v| + 1)}\left(\frac{d(v)(d(v)-1)}{2} - \frac{(|S_v|-1)(|S_v|-2)}{2}\right)\\
     &=\frac{d(v)^2-d(v) - |S_v|^2+3|S_v|-2}{2(d(v)-|S_v|+1)}\\
     &=\frac{(d(v)-|S_v|+1)(|S_v|+d(v)-2)}{2(d(v)-|S_v|+1)} =\frac{|S_v|+d(v)-2}{2}
\end{align*}
Thus, from the above analysis, we get that when $s-1 =2$
\begin{equation}\label{eqn20}
    \frac{1}{d(v) - |S_v| + 1} \left( \binom{d(v)}{s-1} - \binom{|S_v| - 1}{s-1} \right) 
    = \frac{1}{|S_v| + d(v) - 1} \binom{|S_v| + d(v) - 1}{s - 1}
\end{equation}
Thus from \cref{eqn19,eqn20} we get that for $d(v) \geq |S_v| \geq 2$ and $s \geq 3$;
\begin{equation}\label{new19}
    \frac{1}{(d(v)- |S_v| + 1)} \left( {d(v) \choose s-1} - {|S_v|-1 \choose s-1}  \right) \leq \frac{1}{(|S_v|+d(v) - 1)}{|S_v| + d(v) -1\choose s-1}
\end{equation}
Now, adding the inequalities in \cref{y=s ineq} and \cref{new19}, and using them along with \cref{LHS,RHS}, we get that; 
\begin{align*}
    LHS &=  \frac{1}{d(v) - |S_v| + 1} \left( {d(v)  \choose s} - {|S_v|-1 \choose s}+{d(v) \choose s-1} - {|S_v|-1 \choose s-1}\right )\\
    &\leq \frac{1}{{|S_v| + d(v) - 1}} \left({ |S_v| + d(v)-1 \choose s } + { |S_v| + d(v)-1 \choose s -1}\right) = RHS
\end{align*}
Thus we conclude that the bound in \cref{enough2} holds when $|S_v| > 1$.\end{case}

\vspace{2mm}
\noindent We addressed the validity of \cref{thm1} for $s = 1$ in \cref{s=1}, and assumed $s \geq 2$. 
At this point, it is convenient to distinguish between the cases $s = 2$ and $s \geq 3$.

\paragraph{The case $s=2$.} 
It is straightforward to verify that the bound in \cref{enough2} holds with equality when $s = 2$. 
However, the inequalities that arise in the proof of \cref{enough2} for $s \geq 3$ provide additional leverage in ruling out certain cases during the analysis of extremal graphs. 
This advantage is absent when $s = 2$. 
Thus, although the required inequality is easy to establish for $s=2$, the identification of the corresponding extremal class requires a substantially more delicate analysis. 
Consequently, our method, while sufficient to prove the inequality, does not contribute to the extremal characterization in this case. 
For a complete treatment of the $s=2$ case, we refer the reader to \cite{adak2025vertex}.

\paragraph{The case $s \geq 3$.} 
For $s \geq 3$, the bound in \cref{enough2} holds whenever $|S_v| \geq 1$, thereby establishing \cref{main_claim} in this range. 
The remaining case of $|S_v| = 0$ was already handled in \cref{|S_v|=0}.

\medskip
In summary, \cref{main_claim} is established in full generality for all $s \geq 2$.

\end{proof}
\end{lem}
\noindent We proceed according to the following iterative algorithm:
\begin{enumerate}
    \item Initialize $G_0 = G$ and $L_0 = L(G_0, P, u)$.
    
    \item Define $G_1' = G_0 \setminus L_0$, and remove all isolated vertices from $G_1'$ to obtain $G_1$.
    
    \item Set the iteration index $i = 1$ and initialize $x = u$.
    
    \item \textbf{While} $V(G_i) \neq \emptyset$, perform the following steps:
    \begin{enumerate}
        \item For each $v \in V(G_i)$, let $c_i(v)$ denote the weight of vertex $v$ in the context of $G_i$.
        
        \item \textbf{If} $x \notin V(G_i)$:
        \begin{enumerate}
            \item Update $x$ to be a vertex in $G_i$ with maximum weight, that is, $c_i(x) \geq c_i(v)$ for all $v \in V(G_i)$
        \end{enumerate}
        
        \item Let $P_i$ be a longest $x$-path in $G_i$.
        
        \item Define $L_i = L(G_i, P_i, x)$.
        
        \item Define $G_{i+1}' = G_i \setminus L_i$.
        
        \item Remove all isolated vertices from $G_{i+1}'$ to obtain $G_{i+1}$.
        
        \item Increment the iteration index: $i \leftarrow i+1$.
    \end{enumerate}
\end{enumerate}
Let the algorithm run for $t$ iterations, where $t = \max\{i \mid V(G_i) \neq \emptyset\}$ denotes the final iteration index for which the graph $G_i$ remains non-empty.
\begin{rem}\label{u notin L}
    Let $j = \max\{i \mid u \in V(G_i)\}$. Then, by construction, $u \notin L_k$ for all $k \leq j-1$, since $u$ is still present in the graph $G_j$ at step $j$. Moreover, as the path $P_j$ is chosen to start at $u$, and $L_j$ is defined as the set of terminal vertices of transforms of $P_j$, excluding $u$, it follows that $u \notin L_j$. Since $u \notin V(G_{j+1})$, it is removed at step $j$ (as an isolated vertex of $G_{j+1}'$), and hence cannot belong to $L_k$ for any $k \geq j+1$. Therefore, we conclude that $u \notin L_i$ for all $i \in [t]\cup\{0\}$.
\end{rem}
\begin{obs}
For all $i \in [t]$, a bound analogous to \cref{main_claim} holds for $G_i$, with parameters $c_i(v)$ replacing $c(v)$ and $L_i$ replacing $L$. Specifically, we have:
\begin{equation}\label{eqn21}
    N(G_i, K_s, L_i) \leq \sum_{v \in L_i} \frac{1}{c_i(v) - 1} {c_i(v) \choose s}
\end{equation}
\end{obs}

\noindent It is easy to verify that;
\begin{equation}\label{eqn22}
    N(G,K_s) = \sum_{i=0}^tN(G_i,K_s,L_i)
\end{equation}
Note that for all $i\in [t]$ and $v \in V(G_i)$, \begin{equation}\label{ci<c}
    c_i(v) \leq c(v)
\end{equation}
Now from \cref{eqn21,eqn22} and using \cref{inc} we get that;
\begin{equation}\label{algo}
    N(G,K_s) \leq \sum_{i=0}^t\sum_{v \in L_i}\frac{1}{c_i(v)-1}{c_i(v) \choose s} \leq \sum_{i=0}^t\sum_{v \in L_i}\frac{1}{c(v)-1}{c(v) \choose s}
\end{equation}
Since $u \notin L_i$ for any $i \in [t]\cup\{0\}$, we have $\bigcupdot_{i=0}^tL_i \subseteq V(G)\setminus \{u\}$. Therefore from \cref{algo} we get;
\begin{equation}\label{eqn25}
    N(G,K_s) \leq \sum_{i=0}^t\sum_{v \in L_i}\frac{1}{c(v)-1}{c(v) \choose s} \leq \sum_{V(G)\setminus \{u\}}\frac{1}{c(v)-1}{c(v) \choose s}
\end{equation}

\subsubsection{Preprocessing for characterization of Extremal Graphs}
 Recall that the set $H_c$ consists of all vertices in $G$ whose weights are at least $s$.
\begin{lem}
The graph $G$ is extremal for \Cref{thm1} if and only if the induced subgraph $G[H_c]$ is extremal. 
\begin{proof}This follows from the fact that vertices in the complement set $H_c^c = V(G) \setminus H_c$ do not contribute to either side of the inequality in \Cref{thm1}. Note that $N(G,K_s) = N(G[H_c],K_s)$, since any clique of order at least $s$ in $G$, that is not entirely in $G[H_c]$ has to contain a vertex $v \in H_c^c$; but then $c(v) \geq s$, contradiction. On the other hand, since for $v \in H_c^c$, we have $c(v) < s$, we get ${c(v) \choose s} = 0$. 
\end{proof}
\end{lem}
 \begin{rem}\label{heavyG}
     To characterize the extremal graph, we will assume that $G$ does not contain any vertex with weight less than $s$.
 \end{rem}

 \begin{notn}
     Let $v \in V(L_i)$, then the set of neighbors of $v$ in $G_i$, which are outside $L_i$ is denoted by $S^i_v$ and let $d_i(v)$ denote the degree of the vertex $v$ in $G_i$.
 \end{notn}

\begin{lem}\label{s0}
    Let $G$ be an extremal graph for \Cref{thm1}. Then, for every $i \in [t] \cup \{0\}$, there does not exist any vertex $v_0 \in L_i$ such that $|S^i_{v_0}| = 0$.

    \begin{proof}
        Suppose $v_0 \in L_i$ with $|S^i_{v_0}| = 0$ for some $i$. Since $G$ is extremal, the bound in \cref{sv =0} must be an equality, that is
        \begin{equation*}
            \frac{1}{c_i(v_0)}{c_i(v_0) \choose s} = \frac{1}{c_i(v_0) -1}{c_i(v_0) \choose s}
        \end{equation*}
        Thus, both sides must be $0$, that is, $c_i(v_o)<s$. Also, we must have equality in \cref{algo}, and thus in \cref{ci<c} as well. Therefore, $c(v_0) = c_i(v_0) <s$. Thus from \cref{heavyG} we get, $v_0 \notin V(G)$.
    \end{proof}
\end{lem}
\begin{obs}
    Since $G$ is extremal, we must have equality in \cref{s<c-l,dv<l} for the graph $G_i$ with parameters $|L_i|, d_i(v)$ and $|S^i_v|$; that is for all $v \in L_i$,
\begin{equation}\label{boundsLi}
    |S^i_v|+|L_i| = c_i(v) \text{ and } d_i(v) = |L_i|
\end{equation}
\end{obs}
\begin{lem}\label{sg1}

Let $G$ be an extremal graph for \Cref{thm1}. Then, for every $i \in [t] \cup \{0\}$, there does not exist any vertex $v_0 \in L_i$ such that $|S^i_{v_0}| >1$.

    \begin{proof}
        Suppose $v_0 \in L_i$ such that $|S^i_{v_0}| >1$ for some $i$. The bound in \cref{y=s ineq} must be tight. Therefore, from \cref{y=s equality} and using \cref{boundsLi} we get that $|S^i_{v_0}| + d_i(v_0) -1 = |S^i_{v_0}| + |L_i| -1 = c_i(v_0) -1 < s$.
        
        \vspace{1mm}
        \noindent Recall that $s \geq 3$. Now suppose $s-1 =2$. Therefore, $c_i(v_0) \leq s = 3 \implies |S^i_{v_0}| + d_i(v_0) \leq 3$. Since $|S^i_{v_0}| > 1 \implies d_i(v_0)\leq 1$. But $|S^i_{v_0}| \leq d_i(v_0)$. Thus we get a contradiction. Therefore $s-1 \geq 3$.

        \vspace{1mm}
        \noindent Since $G$ is extremal we must have equality in \cref{eqn19}. Thus from \cref{y=s-1 equality} we get, $|S^i_{v_0}|+d_i(v_0)-1 < s-1 \implies |S^i_{v_0}| + d_i(v_0) <s$. Since, $|S^i_{v_0}| + d_i(v_0) = |S^i_{v_0}| + |L_i| = c_i(v_0)$, from equality in \cref{ci<c}, we get that $c(v_0) = c_i(v_0) < s$. Thus from \cref{heavyG} we get that, $v_0 \notin V(G)$.
\end{proof}
\end{lem}
\begin{obs}
Recall that when $s = 2$, the bound in \cref{enough2} always holds with equality if there exists a vertex $v_0 \in L_i$ such that $|S_{v_0}^i| > 1$. Consequently, one cannot derive a contradiction to the existence of such a vertex $v_0$, as was done in \cref{sg1}. To overcome this, in~\cite{adak2025vertex}, the authors had to adopt a more sophisticated structural approach to show that no such vertex exists.
\end{obs}
\subsubsection{Characterization of Extremal graphs}
\vspace{2mm}
\noindent\textbf{$\bullet$ If part: If $G$ is Parent-dominated then $G$ is extremal}
\vspace{2mm}

\noindent We begin by assuming that $G$ is a parent-dominated block graph and aim to show that equality holds in \cref{thm1} for the graph $G$. Let $\{B_1, B_2, \dots, B_k\}$ denote the set of blocks in $G$, where $|V(B_i)| = b_i$ for each $i \in [k]$. Without loss of generality, assume that $B_k$ is the root block; then $B_k$ is a block with the largest order in $G$.

\vspace{1mm}
\noindent Observe that for any vertex $v$, we have $c(v) = \max\{ |B_i| \mid v \in B_i \}$. In particular, all vertices in the root block $B_k$ attain the maximum weight. For each $i \in [k-1]$, let $w_i$ denote the cut-vertex that connects the block $B_i$ to its parent block. Let $u \in V(B_k)$, clearly it is a vertex of maximum weight in $G$. Define $w_k = u$. Since $G$ is a parent-dominated block graph, we have $c(v) = b_i$ for every vertex $v \in V(B_i) \setminus \{w_i\}$. Also note that
\begin{equation}\label{eqn26}
V(G)\setminus \{u\} = \bigcup_{i=1}^k \left( V(B_i) \setminus \{w_i\} \right).
\end{equation}
Recall that the case $s=1$ was addressed in \cref{s=1}, thus assuming $s>1$ we get that;
\begin{equation}\label{eqn27}
    N(G, K_s) = \sum_{i = 1}^k N(B_i, K_s)
\end{equation}
since each $K_s$ lies entirely within a block.

\vspace{1mm}
\noindent We proceed by counting the number of copies of $K_s$ in $B_i$ contributing to the left side of the bound of \cref{thm1}, simultaneously keeping track of the contribution of the weights of the vertices within $B_i$ to the right side of the bound, excluding the vertex $w_i$. Since each $B_i$ is a clique of order at least $s$ (from the assumption in \cref{heavyG}), we have;
\begin{equation}\label{eqn28}
    N(B_i,K_s) = {|V(B_i)| \choose s} = {b_i\choose s}
\end{equation}
Now since $c(v) = b_i$ for all $v \in V(B_i)\setminus\{w_i\}$ we get;
\begin{equation}\label{eqn29}
    \sum_{v \in V(B_i)\setminus\{w_i\}} \frac{1}{c(v)-1}{c(v) \choose s}= \sum_{v \in V(B_i)\setminus\{w_i\}} \frac{1}{b_i-1}{b_i \choose s} = {b_i\choose s} = N(B_i,K_s)
\end{equation}
\noindent Now, from \cref{eqn26,eqn27,eqn28,eqn29}, we obtain:
\begin{equation*}
    \sum_{v \in V(G)\setminus\{u\}} \frac{1}{c(v)-1}{c(v) \choose s}=\sum_{i=1}^k\left(\sum_{v \in V(B_i)\setminus\{w_i\}} \frac{1}{c(v)-1}{c(v) \choose s}\right)=\sum_{i=1}^kN(B_i,K_s) = N(G,K_s)
\end{equation*}
Thus, when 
$G$ is a parent-dominated block graph, equality holds in \Cref{thm1}.

\vspace{2mm}

\noindent\textbf{$\bullet$ Only if part: If $G$ is extremal then $G$ is Parent Dominated}
\vspace{2mm}

    \noindent From now onwards we will assume that $G$ is extremal for \Cref{thm1}. We still follow the assumption of \cref{heavyG}.
\begin{rem}
    Since $G$ is an extremal graph, the second inequality of \cref{eqn25} must be an equality which results $\bigcup_{i=0}^tL_i = V(G) \setminus \{u\}$.  from \cref{s0,sg1} it is clear that for all $i \in [t]\cup\{0\}$ and for all $v \in L_i, |S^i_v| =1$.
\end{rem}
\begin{lem}
    $c_i(v)$ is same for all $v \in L_i$.
    \begin{proof}
        Since $G$ is extremal, from \cref{boundsLi} we have, for all $v \in L_i$, $d_i(v) = |L_i|$. Now since $|S^i_v| = 1$ for all $v \in L_i$, we get that $c_i(v) = d_i(v) +1 = |L_i| +1$.
    \end{proof}
\end{lem}
\begin{rem}
From \cref{No_early_nbr}, it follows that only the vertices occurring \emph{after} the pivot vertex are permuted, while the relative positions of the remaining vertices remain fixed. Moreover, by the above lemma, we have that $c_i(v)$ is the same for all $v \in L_i$. Consequently, every transform of $P_i$ admits the same pivot vertex; denote it by $w_i$.
\end{rem}

\begin{obs}
    Note that $w_i \notin L_i$, since $Pred_{w_i}(1)$ on $P_i$, can not be a neighbor of any vertex in $L_i$.
\end{obs}
\begin{lem}\label{clique}
    $L_i\cup \{w_i\}$ induces a clique in $G_i$ 
    \begin{proof}
        Recall that from \cref{boundsLi}, $d_i(v) = c_i(v)-|S^i_v|$. Since $|S_v^i| =1$, we have, $d_i(v) = c_i(v)-1$. Therefore, $L_i = Back(P_i)$. Thus for all $v \in L_i$, $N(v)\cap V(G_i) = Back^*(P_i)\setminus \{v\}$. Thus $L_i\cup\{w_i\}$ induces a clique in $G_i$. 
    \end{proof}
\end{lem}

\begin{lem}
    If $G$ is an extremal graph for \Cref{thm1}, then $G$ is a parent-dominated block graph.
    \begin{proof}

Proof by induction on the number of vertices.

\noindent\textit{Base Case:} If $|V(G)| = 1$, then $G$ is trivially a parent-dominated block graph.

\noindent\textit{Induction Hypothesis:} Assume that any extremal graph for \Cref{thm1} with fewer than $|V(G)|$ vertices is a parent-dominated block graph.

\noindent\textit{Induction Step:} Let $G$ be an extremal graph on $n = |V(G)|$ vertices. If $G$ is a clique, then it is trivially a parent-dominated block graph and we are done. So, assume $G$ is not a clique.

\vspace{1mm}
\noindent Suppose, for contradiction, that $G$ is disconnected. Without loss of generality, assume $G$ has two components. Let $G'$ and $G''$ be the two components of $G$. Since $G$ is extremal, both $G'$ and $G''$ must be extremal for \Cref{thm1} applied to $G'$ and $G''$ respectively.
\newline Let $u'$ and $u''$ be the heaviest vertices of $G'$ and $G''$, respectively, and clearly $c(u) = \max\{c(u'), c(u'')\}$. Then we have:
\[
N(G',K_s) = \sum_{v \in V(G')\setminus\{u'\}} \frac{1}{c(v)-1} \binom{c(v)}{s}
\]
\[
N(G'',K_s) = \sum_{v \in V(G'')\setminus\{u''\}} \frac{1}{c(v)-1} \binom{c(v)}{s}
\]

Adding these gives:
\[
N(G, K_s) = N(G', K_s) + N(G'', K_s) = \sum_{v \in V(G)\setminus\{u', u''\}} \frac{1}{c(v)-1} \binom{c(v)}{s}
\]

Since $c(u'),c(u'') \geq s$, the above quantity is strictly less than:
\[
\sum_{v \in V(G)\setminus\{u\}} \frac{1}{c(v)-1} \binom{c(v)}{s}
\]
which contradicts the extremality of $G$. Hence, $G$ must be connected.

\vspace{1mm}

\noindent Since $G$ is extremal, each of the $G_i$ must be extremal for \cref{thm1} applied to $G_i$. Therefore, by the induction hypothesis, the graph $G_1 = G \setminus L_0$ is a parent-dominated block graph. We know that $L_0 \cup \{w_0\}$ induces a clique in $G$ by \cref{clique}, where $w_0$ is the only vertex in $G_1$ that is adjacent to the vertices in $L_0$. Thus, $G$ is a block graph, with $B = Back^*(P)$ as a leaf block and ${w}$ as its cut vertex. Let $B'$ be the parent block of $B$ in $G$. Thus $V(B) \cap V(B') = \{w_0\}$. From equality in \cref{ci<c}, the weight of $w_0$ should remain unchanged upon the removal of $L_0$, that is, $c(w_0) = c_1(w_0)$. Since $G_1$ is a parent-dominated block graph, $c(w_0) = |V(B')|$. As the vertices in $B$ induce a clique and $w_0 \in V(B)$, we get $c(w_0) \geq |V(B)|$. It follows that $|V(B')| \geq |V(B)|$. This implies that $G$ is a parent-dominated block graph.
\end{proof}
\end{lem}
\begin{obs}
Note that we characterized the extremal graph under the assumption of \cref{heavyG}. 
A natural question that arises is: what does an extremal graph for \cref{thm1} look like 
when the assumption of \cref{heavyG} is not imposed? 
Let $C$ be a component of $G[H_c^c]$. Note that $C$ may be any graph whose vertices all have weight less than $s$.
It is straightforward to observe that $\lvert N(C) \cap H_c \rvert \leq 1$, 
where $N(C)$ denotes the set of neighbors of the vertices in $C$. 
Indeed, if $C$ had two or more neighbors in $H_c$, then some vertices of $C$ would lie on a cycle of length at least $s$, 
contradicting the fact that their weight is less than $s$. 

\end{obs}
\subsection{Proof of \Cref{thm2}}
We follow a similar inductive strategy as in the proof of \Cref{th:pathl} from \cite{adak2025vertex}.
\subsubsection{Preprocessing steps}
Recall that $H_p = \{v \in V(G) \mid p(v) \geq s-1\}$.
\begin{lem}
The graph $G$ is extremal for \Cref{thm2} if and only if the induced subgraph $G[H_p]$ is extremal.
\begin{proof}This follows from the fact that vertices in the complement set $H_p^c = V(G) \setminus H_p$ do not contribute to either side of the inequality in \Cref{thm2}. Note that $N(G,K_s) = N(G[H_p],K_s)$, since any clique of order at least $s$ in $G$, that is not entirely in $G[H_p]$ has to contain a vertex $v \in H_p^c$; but then $p(v) \geq s-1$, contradiction. On the other hand, since for $v \in H_p^c$, we have $p(v) < s-1$, we get ${p(v) \choose s-1} = 0$. 
 
\end{proof}
\end{lem}
\begin{rem}\label{heavyGpath}
    To characterize the extremal graph, we will assume that $G$ does not contain any
vertex with weight less than $s-1$.
\end{rem}
\begin{rem}
    When $s =1$, the quantity $N(G,K_s)$ simply counts the number of vertices in $G$, therefore $N(G,K_s) = n$. In this case, the right-hand side of the bound in \Cref{thm2} simplifies to:
    \begin{equation*}
        \sum_{v\in V(G)}{p(v) \choose 0} = \sum_{v\in V(G)}1 = n
    \end{equation*}
    Thus, the bound holds with equality for all $G$.
\end{rem}
\begin{rem}
    For $s=2$, \Cref{thm2} coincides with \Cref{th:pathl}. Indeed, since $N(G,K_2)=m$, the number of edges in $G$, the bound in \Cref{thm2} reduces to
    \begin{equation*}
        m \leq \frac{1}{2}\sum_{v \in V(G)} p(v).
    \end{equation*}
    Note that $p(v)=0$ for isolated vertices, so $H_p$ consists precisely of the non-isolated vertices of $G$. For the case $s=2$, we therefore use to the proof of \Cref{th:pathl} from~\cite{adak2025vertex}, which shows that the extremal graphs are exactly the disjoint unions of cliques.
\end{rem}
In what follows, we assume $s >2$.
\subsubsection{Inductive Proof}
\begin{proof} We will prove \Cref{thm2} using induction on the number of vertices in graph $G$.
    
    \vspace{1mm}
    \noindent\textit{Base Case:} If $n = 1$, clearly both sides are zero.

    \vspace{1mm}
    \noindent\textit{Induction Hypothesis:} Assume that the claim holds for all graphs with at most $|V(G)| -1$ vertices, where $|V(G)| \geq 2$.

     \vspace{1mm}
    \noindent\textit{Induction Step:} If $G$ is disconnected by induction hypothesis the statement of the theorem is true for each connected component and thus for $G$. Therefore assume $G$ to be connected. Let $k = max\{p(v) \mid v \in V(G)\}$ and $P$ be a $k$-length path in $G$. Let $P = v_0,v_1,\dots,v_k$. Now consider the following two cases;
    \begin{case}\label{case1}
        There exists a $(k+1)$-length cycle in $G$: We will show that in this case the extremal graph has to be a clique.
        \newline Without loss of generality assume the endpoints of $P$, $v_0$ and $v_k$, are adjacent thus forming a $(k+1)$-length cycle.
        \begin{claim}
     $V(P) = V(G)$.
     \begin{proof}
          We know that $v_0$ is adjacent to $v_k$. Suppose $V(G) \setminus V(P) \neq \emptyset$. Since $G$ is connected, there exists a vertex  $v \in V(G)\setminus V(P)$, adjacent to some vertex in $P$. Since the vertices of $P$ form a cycle, without loss of generality, we can assume $v$ is adjacent to $v_0$. Now consider the path;
         \begin{equation*}
             Q = v,v_0,v_1,\dots,v_{k}
         \end{equation*}
Thus $Q$ is a path in $G$ with length $k +1$. Thus we arrive at a contradiction to the assumption that the length of a longest path in $G$ is $k$.
     \end{proof}
 \end{claim}
 \noindent From the above claim we get that $p(v) = |V(G)| -1 =(n-1)$ for all $v \in V(G)$. Therefore;
 \begin{equation}\label{eqn1}
     \sum_{v \in V(G)}\frac{1}{s}{p(v) \choose s-1} = \frac{n}{s}{n-1 \choose s-1} = {n \choose s}
 \end{equation} 
 \noindent Note that the number of copies of $K_s$ in $G$ is at most the number of $s$-subsets of $V(G)$. Thus;
 \begin{equation}\label{eqn2}
     N(G,K_s) \leq {n \choose s}
 \end{equation}
 \noindent From \cref{eqn1} and \cref{eqn2} we get the inequality of \Cref{thm2}. Note that for equality, \cref{eqn2} must be an equality. Thus for all $S \subseteq V(G)$ with $|S| = s$, we need $G[S] \cong K_s$. Since $s >1$, $G$ must be a clique.
    \end{case}
    \begin{case}\label{pathcase2}
        There does not exist any $(k+1)$-length cycle in $G$: We will show that in this case, no extremal graph exists.
        \newline From the assumption of the case, $v_0$ and $v_k$ are non-adjacent. Without loss of generality assume $d(v_0) \geq d(v_k)$.
\begin{claim}\label{light}
    $d(v_k) \leq \frac{k}{2}$
    \begin{proof}
        Suppose not, that is $d(v_k) \geq \frac{k+1}{2}$. Thus we get $d(v_0) \geq \frac{k+1}{2}$ also. Since $P$ is a longest path in $G$ and $v_0$ is not adjacent to $v_{k}$, we get, $N(v_0), N(v_k) \subseteq \{v_1,v_2,\dots, v_{k-1}\}$. Let $N = N(v_k) \setminus \{v_{k-1}\}$ and $N^+ = \{v_{j+1} \mid v_j \in N\}$. Clearly, $|N| = |N^+| \geq \frac{k+1}{2}-1 = \frac{k-1}{2}$. Note that $v_0, v_k \notin N^+$. Therefore, $|\{v_1,v_2,\dots,v_{k-1}\} \setminus N^+| \leq \frac{k-1}{2}$. Since $|N(v_0)| \geq \frac{k+1}{2}$, by pigeonhole principle, we get that, there exists $v_i \in N^+$ such that $v_i \in N(v_0)$ and $v_{i-1} \in N(v_k)$. Consider the cycle; $$K = v_0, v_1, \dots, v_{i-1},v_k,v_{k-1},\dots, v_i,v_0$$ Clearly $K$ is a $(k+1)$-length cycle in $G$. Thus, we get a contradiction to the assumption of \cref{pathcase2}.
    \end{proof}
    \end{claim}
\begin{notn}
    Let $\mathcal{K} = \lfloor\frac{k}{2}\rfloor$.
\end{notn}
\begin{rem}\label{lighter}
Thus from \cref{light} and since $d(v_k) \in \mathbb{Z}$, we get $d(v_k) \leq \mathcal{K}$.
\end{rem}
\noindent Let $G' = G \setminus \{v_k\}$. Let $p_{G'}(v)$ denote the weights of the vertices restricted to the graph $G'$. Thus using the induction hypothesis and the fact that $p_{G'}(v) \leq p(v)$ for all $v \in V(G')$ we get;
\begin{equation}\label{eqn5}
    N(K_s,G') \leq \sum_{v \in V(G')}\frac{1}{s}{p_{G'}(v) \choose s-1} \leq \sum_{v \in V(G')}\frac{1}{s}{p(v) \choose s-1}
\end{equation}
Recall that $N(G,K_s,\{v\})$ denotes the number of copies of $K_s$ in $G$, containing the vertex $v$. Clearly;
\begin{equation}\label{eqn4}
    N(G,K_s) = N(G',K_s) + N(G,K_s,\{v_k\})
\end{equation}
Therefore, to show that the bound in \Cref{thm2} is true, using \cref{eqn5,eqn4}, it is enough to show that;
\begin{equation}\label{eqn6}
    N(G,K_s,\{v_k\}) \leq \frac{1}{s}{p(v_k) \choose s-1}
\end{equation}
\begin{rem}If $N(G,K_s,\{v_k\}) = 0$, clearly \cref{eqn6} holds and we are done. Thus assume $N(G,K_s,\{v_k\}) > 0$. Therefore, $d(v_k) \geq (s-1)$. Thus from \cref{lighter},  we have $\mathcal{K}\geq (s-1)$.
\end{rem}
\noindent Note that for every $K_s$ containing $v_k$, there exists a unique copy of $K_{s-1}$ in $G[N(v_k)]$. Therefore;
\begin{equation}\label{eq3}
  N(G,K_s,\{v_k\}) = N(G[N(v_k)],K_{s-1}) \leq {d(v_k) \choose s-1}
\end{equation}

    \noindent From \cref{lighter} and \cref{eq3} we get;
    \begin{equation}\label{eqn7}
         N(G,K_s,\{v_k\}) \leq {d(v_k) \choose s-1}\leq {\mathcal{K}\choose s-1}
    \end{equation}
    \begin{obs}From \cref{eqn6,eqn7}, for the bound of \Cref{thm2} to be true, it is enough to show that;
    \begin{align}
        &{\mathcal{K}\choose s-1} \leq \frac{1}{s}{p(v_k) \choose s-1} = \frac{1}{s}{k \choose s-1}\nonumber\\
        \iff& \frac{\left(\mathcal{K}\right)!}{(s-1)!\left(\mathcal{K}-s+1\right)!} \leq \frac{k!}{s(s-1)!(k-s+1)!}\nonumber\\
        \iff& s\left(\mathcal{K}\right)\left(\mathcal{K}-1\right)\dots\left(\mathcal{K}-s+2\right) \leq k(k-1)\dots(k-s+2)\nonumber\\
        \iff&s \leq \frac{k}{\left(\mathcal{K}\right)}\frac{k-1}{\left(\mathcal{K}-1\right)}\dots\frac{k-s+2}{\left(\mathcal{K}-s+2\right)}\label{eqn8}
    \end{align}
    \end{obs}
    \begin{rem}\label{>2}
        Since $\mathcal{K} = \lfloor\frac{k}{2}\rfloor$, for any $x > 0$ and $\mathcal{K} > x$ we have;
        \begin{equation*}
            \frac{k-x}{\mathcal{K}-x} > \frac{k}{\mathcal{K}} \geq 2    
        \end{equation*}
    \end{rem}

\noindent From \cref{>2} and the fact that $s>2$ we get;
\begin{equation}\label{*}
    \frac{k}{\left(\mathcal{K}\right)}\frac{k-1}{\left(\mathcal{K}-1\right)}\dots\frac{k-s+2}{\left(\mathcal{K}-s+2\right)} > 2^{s-1}
\end{equation}
It is easy to see that $2^{s-1} > s$ for $s >2$. Thus from \cref{*} we get that \cref{eqn8} holds;
\begin{equation*}
    s < 2^{s-1} < \frac{k}{\left(\mathcal{K}\right)}\frac{k-1}{\left(\mathcal{K}-1\right)}\dots\frac{k-s+2}{\left(\mathcal{K}-s+2\right)}
\end{equation*}
Therefore, the bound in \Cref{thm2} is true, but equality does not hold in this case.

\end{case}
\noindent Thus, the bound for \cref{thm2} holds for all $G$, but equality holds only for \cref{case1}, that is, equality holds if and only if all the components of $G$ are cliques.
\end{proof}
\begin{obs}
We characterized the extremal graph under the assumption of \Cref{heavyGpath}. Without this assumption, it is easy to see that no vertex in $H_p^c$ can be adjacent to a vertex of $H_p$. Consequently, the components of $G$ are either cliques of order at least $s$, which belong to $H_p$, or arbitrary graphs whose vertex weights are all less than $s-1$.

\end{obs}
\section{An Additional Result}

\begin{notn}
    For a graph $G$, let $l(v)$ denote the weight of a vertex $v \in V(G)$, defined as
    \[
        l(v) = \text{length of the longest $v$-path in $G$}.
    \]
\end{notn}

\begin{thm}
    \emph{(Balister--Bollob\'{a}s--Riordan--Schelp \cite{BALISTER2003366})}\label{balbol} 
    For every graph $G$, we have
    \[
        |E(G)| \leq \sum_{v \in V(G)} \frac{l(v)}{2}.
    \]
    Equality holds if and only if every connected component of $G$ is a clique.
\end{thm}

It is worth noting the close resemblance between \Cref{balbol} and \Cref{th:pathl}. Indeed, since $l(v) \leq p(v)$ for every $v \in V(G)$, Theorem~\ref{th:pathl} follows directly from Theorem~\ref{balbol}.

Moreover, one can replicate the proof of \Cref{thm2} by replacing $p(v)$ with $l(v)$ throughout, and the argument remains valid. This yields the following result, which both generalizes Theorem~\ref{balbol} and strengthens \Cref{thm2}.

\begin{thm}
    Let $G$ be a simple graph. Then,
    \[
        N(G,K_s) \leq \sum_{v \in V(G)} \frac{1}{l(v)+1} {l(v)+1 \choose s} 
        = \frac{1}{s} \sum_{v \in V(G)} {l(v) \choose s-1}.
    \]
    If $s = 1$, equality always holds. Otherwise, equality holds if and only if every component of $G[H_l]$ is a clique, where
    \[
        H_l = \{v \in V(G) \mid l(v) \geq s-1\}.
    \]
\end{thm}

\bibliographystyle{plain}
\bibliography{references}
\end{document}